\documentclass[11pt,twoside]{article}
\usepackage{amsfonts}
\usepackage{fancyhdr}
\usepackage{titlesec}
\usepackage{cite}
\usepackage{ifthen}
\usepackage{amssymb}
\usepackage{fancyhdr}
\usepackage{titlesec}
\usepackage{pifont}
\usepackage{setspace}
\usepackage{indentfirst}
\usepackage{amsmath,amssymb,amscd,amsthm,mathrsfs}
\input amssym.def

\newboolean{first}
\setboolean{first}{true}
\renewcommand{\headrulewidth}{0pt}

\newfont{\aaa}{cmb10 at 19pt}
\newfont{\bbb}{cmb10 at 11pt}
\newtheorem{lemma}{Lemma}[section]
\newtheorem{theorem}{Theorem}[section]
\newtheorem{definition}{Definition}[section]

\newcommand{\beq}{\begin{equation}}
\newcommand{\eeq}{\end{equation}}
\newcommand{\bey}{\begin{eqnarray}}
\newcommand{\eey}{\end{eqnarray}}
\newcommand{\beyy}{\begin{eqnarray*}}
\newcommand{\eeyy}{\end{eqnarray*}}

\numberwithin{equation}{section}
\makeatletter
\def\@oddfoot{}
\makeatother

\begin{document}
\thispagestyle{empty} \thispagestyle{fancy} {
\fancyhead[RO,LE]{\scriptsize \bf 
} \fancyfoot[CE,CO]{}}
\renewcommand{\headrulewidth}{0pt}

\begin{center}
{\bf \LARGE A Generalization of Stampacchia Lemma\\
 and Applications}

 \vspace{5mm}

{\small \textsc{Hongya GAO} \quad \textsc{Jiaxiang ZHANG}  \quad \textsc{Hongyan MA}\footnote{Corresponding author, email:
mahongyan@hbu.edu.cn.}

\vspace{2mm}

{\small College of Mathematics and Information Science, Hebei University, Baoding, 071002, China
}}
\end{center}

\vspace{3mm}

\begin{center}
\begin{minipage}{135mm}
{\bf \small Abstract}.  {\small We present a generalization of Stampacchia Lemma and give applications to regularity property of weak and entropy solutions of degenerate elliptic equations of the form
$$
\left\{
\begin{array}{llll}
-\mbox{div} (a(x,u(x)) Du (x)) =f(x), & \mbox { in } \Omega, \\
u(x)=0,  & \mbox { on } \partial \Omega,
\end{array}
\right.
$$
where
$$
\frac {\alpha}{(1+|u|) ^\theta} \le a(x,s)\le \beta
$$
with $0<\alpha \le \beta <\infty$ and $0\le \theta <1$.}


\end{minipage}
\end{center}

\thispagestyle{fancyplain} \fancyhead{}
\fancyhead[L]{\textit{}\\
} \fancyfoot{} \vskip 6mm

\section{Stampacchia Lemma and a generalization.}
We first recall the classical Stampacchia Lemma, which can be found, for example, in Lemma 4.1 of \cite{St2}.

\begin{lemma}\label{Stampacchia Lemma}
Let $c, \alpha, \beta$ be positive constants and $k_0\in \mathbb R$. Let
$\varphi: [k_0,+\infty)$ $\rightarrow [0, +\infty)$ be nonincreasing
and such that
\begin{equation}\label{stampacchia_assumption}
\varphi (h) \le \frac {   c}{(h-k)^\alpha} [\varphi (k)]^\beta
\end{equation}
\noindent
for every $h,k$ with $h>k\ge k_0$.
It results that:

 {\bf (i)} if $\beta > 1$ then
\begin{equation*}
 \label{statement_stampacchia_beta>1}
\varphi(k_0 + d) = 0,
\end{equation*}
where
\begin{equation*}
 \label{d=}
d^\alpha =   c [\varphi(k_0)]^{\beta - 1} 2^{\frac {\alpha \beta} {\beta
- 1}}.
\end{equation*}

 {\bf(ii)} if $\beta = 1$ then for any $k\ge k_0$,
\begin{equation*} \label{statement_stampacchia_beta=1}
\varphi (k) \le \varphi(k_0) e^{1 - (  c e)^{-\frac{1}{\alpha}} (k - k_0)}.
\end{equation*}

 {\bf(iii)} if $0< \beta < 1$ and $k_0 > 0$ then for any $k\ge k_0$,
\begin{equation*}
 \label{statement_stampacchia_beta<1}
\varphi (k) \le
2^{\frac {\alpha}{(1-\beta)^2 }}
\left\{   c ^\frac {1}{1-\beta} +(2 k_0)^{\frac \alpha
{1-\beta}} \varphi(k_0) \right\}
\left(\frac 1 k \right) ^{\frac \alpha {1-\beta}}.
\end{equation*}
\end{lemma}
Stampacchia Lemma is an efficient tool in dealing with summability properties of elliptic partial differential equations and systems as well as minima of variational integrals, see \cite{Boccardo1}. Till now, this lemma is
used repeatedly by many mathematicians, see, for example, \cite{ABFOT,BC,BG,BGM,CL,GDM,GLT,IF,K,FS1,FS2}.

Due to the importance of Stampacchia Lemma, some generalizations have been given in order to deal with more complicated situations, see \cite{KV,GLW,GDHR,GHR}.
We emphasize that, in \cite{GLW}, Gao, Leonetti and Wang gave some remarks on Stampacchia Lemma, which compared the condition of Stampacchia Lemma and a special case of it; in \cite{GDHR}, the authors presented two
generalizations of Stampacchia Lemma and gave some applications to elliptic systems;
in \cite{KV}, Kovalevsky and Voitovich dealt with an analogous condition, but contains an additional factor on its right hand side; in \cite{GHR}, Gao, Huang and Ren gave an alternative, seems to be more elementary, proof
of the results due to Kovalevsky and Voitovich \cite{KV}.

For convenience of the reader, we list the results due to Kovalevsky and Voitovich \cite{KV}. The following version can be found in \cite{GHR}.

\begin{lemma}\label{lemma 2.3}
Let $c , \alpha, \beta, k_0$ be positive constants and $0\le \theta <1$. Let
$\varphi: [k_0,+\infty)$ $\rightarrow [0, +\infty)$ be nonincreasing
and such that
\begin{equation}
 \label{generalied stampacchia_assumption}
\varphi (h) \le \frac {c k^{  \theta    \alpha}}{(h-k)^\alpha}
[\varphi (k)]^\beta
\end{equation}
\noindent
for every $h,k$ with $h>k\ge k_0>0$.
It results that:

{\bf(i)} if $\beta > 1$ then there exists $k^*>k_0$ such that $\varphi (k^*)=0$;

{\bf (ii)} if $\beta = 1$ then for any $k\ge k_0$,
\begin{equation*} \label{generalied statement_stampacchia_beta=1}
\varphi (k)\le \varphi(k_0) e ^{1-\left( \frac {k-k_0}{\tau} \right) ^{1-\theta  } },
\end{equation*}
where
\begin{equation*}\label{tau11}
\tau =\max \left\{k_0, (ce2^{\theta   \alpha} (1-\theta  ) ^\alpha)^{\frac 1 {(1-\theta  ) \alpha}} \right\};
\end{equation*}

{\bf  (iii)} if $0< \beta < 1$ then for any $k\ge k_0$,
\begin{equation*}
 \label{statement_stampacchia_beta<1}
\varphi (k) \le
2^{\frac {\alpha (1-\theta  )}{(1-\beta)^2 }}
\left\{c^\frac {1}{1-\beta} +(2 k_0)^{\frac {\alpha (1-\theta  )}
{1-\beta}} \varphi(k_0) \right\}
\left(\frac 1 k \right) ^{\frac {\alpha (1-\theta  )} {1-\beta}}.
\end{equation*}
\end{lemma}

Compare (\ref{generalied stampacchia_assumption}) with (\ref{stampacchia_assumption}), we find that (\ref{generalied stampacchia_assumption}) contains an additional factor depending on $k$ on the numerator of its right hand
side, with the power less than the one on the denominator. This case appears naturally in dealing with degenerate second-order elliptic problems considered in \cite{GHR} and fourth-order elliptic equations considered in
\cite{KV}.

We now present a generalization of Lemma \ref{Stampacchia Lemma}, which can be used in the study of degenerate elliptic equations of divergence type.

\begin{lemma}\label{Stampacchia Lemma, generalization}
Let $  c, \alpha, \beta,k_0$ be positive constants and $0\le \theta <1$. Let
$\varphi: [k_0,+\infty) \rightarrow [0, +\infty)$ be nonincreasing
and such that
\begin{equation}\label{stampacchia_assumption 22}
\varphi (h) \le \frac {c h ^{\theta \alpha }}{(h-k)^\alpha}
[\varphi (k)]^\beta
\end{equation}
\noindent
for every $h,k$ with $h>k\ge k_0 >0$.
It results that:

{\bf (i)} if $\beta > 1$ then
$$
\varphi (2L) =0,
 $$
where
\begin{equation}\label{the value of L}
 L= \max \left\{2k_0,  c^{\frac 1 {(1-\theta ) \alpha}}  [\varphi (k_0)] ^{\frac {\beta -1}{ (1-\theta) \alpha}}  2 ^{\frac 1 {(1-\theta)\beta} \left(\beta +\theta +\frac 1 {\beta -1}\right)} \right\} >0.
\end{equation}

{\bf(ii)} if $\beta = 1$ then for any $k\ge k_0$,
\begin{equation*} \label{generalied statement_stampacchia_beta=1}
\varphi (k)\le \varphi(k_0) e ^{1-\left( \frac {k-k_0}{\tau} \right) ^{1-\theta }},
\end{equation*}
where
\begin{equation}\label{tau11}
\tau =\max \left\{k_0, \left( ce 2^{\frac {(2-\theta ) \theta \alpha }{1-\theta }} (1- \theta ) ^\alpha \right)^{\frac 1 {(1- \theta ) \alpha}} \right\};
\end{equation}

{\bf(iii)} if $0< \beta < 1$ then for any $k\ge k_0$,
\begin{equation} \label{statement_stampacchia_beta<1}
\varphi (k) \le 2^{\frac { (1-\theta)\alpha}{(1-\beta)^2 }} \left\{ (c_1 2^{\theta \alpha})^\frac {1}{1-\beta} +(2 k_0)^{\frac {(1-\theta)\alpha} {1-\beta}} \varphi(k_0) \right\} \left(\frac 1 k \right) ^{\frac
{\alpha(1-\theta)} {1-\beta}},
\end{equation}
where
\begin{equation} \label{c4andc5}
c_1=\max \left\{ 4^{(1-\theta) \alpha} c 2^{\theta \alpha}, c_2^{1-\beta} \right\}, \ \ c_2= 2 ^{\frac {(1-\theta )\alpha}{(1-\beta)^2 }} \left[(c2^{\theta \alpha})^{\frac 1 {1-\beta}} + (2k_0) ^{\frac {(1-\theta)\alpha}
{1-\beta}} \varphi (k_0)\right].
\end{equation}
\end{lemma}

We remark that the difference between (\ref{stampacchia_assumption 22}) and (\ref{generalied stampacchia_assumption}) is that $h^{\theta \alpha}$ in place of $k^{\theta \alpha}$. It is obvious that
(\ref{stampacchia_assumption 22}) is weaker than (\ref{generalied stampacchia_assumption}) since $k<h$. If $\theta =0$, then (\ref{stampacchia_assumption 22}) and (\ref{generalied stampacchia_assumption}) coincide with
(\ref{stampacchia_assumption}).

In order to prove Lemma \ref{Stampacchia Lemma, generalization}, we need two preliminary lemmas. The following lemma can be found in Lemma 7.1 in \cite {Giusti}.

\begin{lemma}\label{iteration lemma}
Let $\beta, B,C,x_i$ be such that $\beta >1$, $C>0$, $B>1$, $x_i \ge 0$  and
\begin{equation}\label{iteration inequality}
x_{i+1} \le CB^i x_i^{\beta}, \quad i=0,1,2,\cdots.
\end{equation}
If
$$
x_0 \le C^{-\frac 1 {\beta -1}} B^{-\frac 1 {(\beta -1)^2 }},
$$
then, we have
$$
x_i \le B^{-\frac i \alpha }x_0, \quad i=0,1,2,\cdots,
$$
so that
$$
\lim _{i\rightarrow +\infty} x_i =0.
$$
\end{lemma}

The following lemma comes from Remark 1 and its proof in \cite{GLW}.

\begin{lemma}\label{GLW}
Let $c_3,\tilde \alpha, \beta, k_0$ be positive constants. Let $\varphi : [k_0, +\infty) \rightarrow [0,+\infty)$ be nonincreasing and such that
\begin{equation}\label{the second condition}
\varphi(2k) \le \frac {c_3}{k^ {\tilde \alpha } } [\varphi (k)] ^{\beta}, \ \  0<\beta <1,
\end{equation}
for every $k\ge k_0$. Then for every $h,k$ with $h>k\ge k_0$,
$$
\varphi (h) \le \frac {c_4}{ (h-k)^{\tilde \alpha}} [\varphi (k)] ^\beta,
$$
where
$$
c_4 =\max \left\{4^{\tilde \alpha} c_3 , c_5^{1-\beta} \right\}, \ \ c_5= 2 ^{\frac {\tilde \alpha}{(1-\beta)^2 }} \left[c_3^{\frac 1 {1-\beta}} + (2k_0) ^{\frac {\tilde \alpha} {1-\beta}} \varphi (k_0)\right].
$$
\end{lemma}

\noindent {\bf Proof of Lemma \ref{Stampacchia Lemma, generalization}}.

(i) For $\beta >1$ we fix $L>k_0$ and choose levels
$$
k_i = 2L (1-2^{-i-1}), \ \  i=0,1,2,\cdots.
$$
It is obvious that $k_0<L\le k_i <2L$ and $\{ k_i \}$ be an increasing sequence. We choose in (\ref{stampacchia_assumption 22})
$$
k=k_i,\  h=k_{i+1}, \ x_i = \varphi (k_i) \  \mbox { and  }  \ x_{i+1} = \varphi (k_{i+1}),
$$
and noticing that
$$
h-k = k_{i+1}-k_i=L 2 ^{-i-1},
$$
we have,
$$
x_{i+1} \le \frac {c \left( 2L (1-2^{-i-2}) \right) ^{\theta \alpha}}{ (L2^{-i-1}) ^\alpha} x_i^\beta,\  \ i=0,1,2,\cdots,
$$
from which we derive
$$
x_{i+1} \le \frac {c 2^{(1+\theta ) \alpha } (2^\alpha ) ^i}{ L^{ (1-\theta ) \alpha}} x_i^\beta, \  \ i=0,1,2,\cdots.
$$
Thus (\ref{iteration inequality}) holds true with (we keep in mind that $\beta >1$)
$$
C=\frac {c 2^{(1+\theta ) \alpha } }{ L^{ (1-\theta ) \alpha}} \  \ \mbox { and } \  \ B=2^\alpha.
$$
We get from Lemma \ref{iteration lemma}
\begin{equation}\label{limit 1}
\lim _{i\rightarrow +\infty} x_i =0
\end{equation}
provided that
\begin{equation}\label{condition 1}
x_0 =\varphi (k_0) =\varphi (L) \le \left(\frac {c 2^{(1+\theta ) \alpha } }{ L^{ (1-\theta ) \alpha}}\right) ^{-\frac 1 {\beta -1}} \left(2 ^\alpha \right)  ^{-\frac 1 {(\beta -1)^2 }} .
\end{equation}
Note that (\ref{limit 1}) implies
$$
\varphi (2L) =0.
$$
Let us check condition (\ref{condition 1}) and determine the value of $L$. (\ref{condition 1}) is equivalent to
\begin{equation}\label{limit 2}
\varphi (L) \le c ^{-\frac 1 {\beta -1}} 2 ^{-\frac {\alpha }{\beta -1 } \left( 1+\theta +\frac 1{\beta -1}\right)}  L ^{\frac { (1-\theta ) \alpha }{\beta -1}}.
\end{equation}
In (\ref{stampacchia_assumption 22}) we take $k=k_0$ and $h=L\ge 2k_0$ (which is equivalent to $L-k_0 \ge \frac L 2$) and we have
$$
\varphi (L) \le \frac {c L^{\theta \alpha }}{(L-k_0)^\alpha } [\varphi (k_0)] ^\beta \le \frac {c 2 ^\alpha L^{\theta \alpha }}{L^\alpha } [\varphi (k_0)] ^\beta = \frac {c2^\alpha [\varphi (k_0)] ^\beta }{ L^{(1-\theta )
\alpha}}.
$$
Then condition (\ref{limit 2}) would be satisfied if
\begin{equation}\label{condition no.1}
L\ge 2k_0
\end{equation}
and
$$
\frac {c 2^\alpha [\varphi (k_0)] ^\beta }{ L^{(1-\theta ) \alpha}} \le c^{-\frac 1 {\beta -1}} 2 ^{-\frac {\alpha }{\beta -1 } \left(1+\theta +\frac 1{\beta -1} \right)}  L ^{\frac { (1-\theta ) \alpha }{\beta -1}},
$$
that is
\begin{equation}\label{condition no.2}
 c^{\frac \beta {\beta -1}} [\varphi (k_0)] ^\beta  2 ^{\alpha + \frac {\alpha }{\beta -1 } \left(1+\theta +\frac 1 {\beta -1} \right)}    \le L^{\frac {(1-\theta ) \alpha \beta }{\beta -1}}.
\end{equation}
 (\ref{the value of L}) is a sufficient condition for both (\ref{condition no.1}) and (\ref{condition no.2}).

(ii)  Let  $\beta =1$ and $\tau$ be as in (\ref{tau11}).
For any $s\in \mathbb N ^+$ we let
$$
k_s =k_0+ \tau  s^{\frac 1 {1-\theta }},
$$
then $\{k_s\}$ is an increasing sequence and
$$
k_{s+1}-k_s =\tau \left[ (s+1) ^{\frac 1 {1-\theta } } -s^{\frac 1 {1-\theta }}\right].
$$
We use Taylor's formula in order to get
\begin{equation}\label{estimate 10}
k_{s+1}-k_s =\tau \left[\frac 1 {1- \theta } s^{\frac { \theta}  {1- \theta }} +\frac { \theta } { (1- \theta )^2} \xi ^{\frac {2  \theta  -1}{1-
 \theta  }}\right]\ge  \frac \tau {1- \theta } s^{ \frac { \theta}  {1- \theta }},
\end{equation}
where $\xi$ lies in the open interval $(s,s+1)$. In (\ref{stampacchia_assumption 22}) we take $\beta =1$, $k=k_s$ and $h=k_{s+1}$, we use (\ref{estimate 10}) and we get
\begin{equation}\label{estimate 11}
\varphi (k_{s+1}) \le \frac { c \left[k_0 +\tau (s+1)^{\frac 1 {1- \theta }}\right] ^{ \theta  \alpha} }{ \left(\frac \tau {1- \theta } \right)^\alpha s^{\frac { \theta  \alpha} {1- \theta }}} \varphi (k_s)
\le \frac { c \left[k_0 +\tau (2s)^{\frac 1 {1- \theta }}\right] ^{ \theta  \alpha} }{ \left(\frac \tau {1- \theta } \right)^\alpha s^{\frac { \theta  \alpha} {1- \theta }}} \varphi (k_s).
\end{equation}
(\ref{tau11}) ensures
$$
k_0\le \tau <\tau (2s) ^{\frac 1 {1- \theta }} \ \ \mbox { and } \ \  \frac { c\left( 2^{1+\frac 1 {1-\theta}}\tau \right) ^{ \theta  \alpha}}{\left( \frac {\tau }{1- \theta }\right) ^\alpha}\le \frac 1 e.
$$
From (\ref{estimate 11}) and the above inequalities, one has
\begin{equation*}\label{estimate 13}
\varphi (k_{s+1}) \le \frac { c \left(2\tau (2s)^{\frac 1 {1- \theta }}\right) ^{ \theta  \alpha} }{ \left(\frac \tau {1- \theta } \right)^\alpha s^{\frac { \theta  \alpha} {1- \theta }}} \varphi (k_s)=\frac { c\left(
2^{1+\frac 1 {1-\theta}}\tau \right) ^{ \theta  \alpha}}{\left( \frac {\tau }{1- \theta }\right) ^\alpha} \varphi (k_s) \le  \frac 1 e \varphi (k_s).
\end{equation*}
 By recursion,
$$
\varphi (k_s) \le \frac 1 {e^s } \varphi (k_0),\ \ s \in \mathbb N^+ .
$$
For any $k\ge k_0$, there exists $s\in \mathbb N ^+$ such that
$$
k_0 +\tau (s-1) ^{\frac 1 {1- \theta }} \le k <  k_0 +\tau s ^{\frac 1 {1- \theta }}.
$$
Thus, considering $\varphi (k)$ is nonincreasing, one has
$$
\varphi (k) \le \varphi \left( k_0 +\tau (s-1) ^{\frac 1 {1- \theta }} \right)=\varphi (k_{s-1}) \le e^{1-s} \varphi (k_0)\le \varphi(k_0) e ^{1-\left( \frac {k-k_0}{\tau} \right) ^{1- \theta }},
$$
as desired.

(iii) We take $h=2k$ in (\ref{stampacchia_assumption 22}) and we get for every $k\ge k_0> 0$,
$$
\varphi(2k) \le \frac {c 2 ^{\theta \alpha}}{ k ^{(1-\theta) \alpha }} [\varphi (k)] ^\beta.
$$
Thus (\ref{the second condition}) holds true with $c_3=c 2 ^{\theta \alpha}$ and $\tilde \alpha =(1-\theta) \alpha$.  Lemma \ref{GLW} tells us that for every $h,k$ with $h>k\ge k_0$,
$$
\varphi(h) \le \frac {c_1}{ (h-k) ^{(1-\theta ) \alpha}} [\varphi (k)]^\beta,
$$
where $c_1$ be as in (\ref{c4andc5}).
We use Lemma \ref{Stampacchia Lemma} (iii) and we get (\ref{statement_stampacchia_beta<1}).

\section{Applications.}
In this section we are interested in the following degenerate elliptic problem
\begin{equation}\label{BVP}
\left\{
\begin{array}{llll}
-\mbox{div} (a(x,u(x)) Du (x)) =f(x), & \mbox { in } \Omega; \\
u(x)=0,  & \mbox { on } \partial \Omega.
\end{array}
\right.
\end{equation}
Here $\Omega$ is a bounded, open subset of $\mathbb R^n$, $n>2$, and $a(x,s):\Omega \times \mathbb R \rightarrow \mathbb R$ is a Carath\'eodory function (that is, measurable with respect to $x$ for every $s\in \mathbb R$,
and continuous with respect to $s$ for almost every $x\in \Omega$) satisfying the following conditions:
\begin{equation}\label{conditon for a(x,s)}
\frac {\alpha}{(1+|s|)^\theta} \le a (x,s) \le \beta,
\end{equation}
for some positive number $\theta $ such that
$$
0\le \theta <1,
$$
for almost every $x\in \Omega$ and for every $s\in \mathbb R$, where $\alpha$ and $\beta$ are positive constants. As far as the source $f$ is concerned, we assume
\begin{equation}\label{condition for f}
f\in L_{weak} ^m(\Omega), \  \ m>1,
\end{equation}
where $L_{weak} ^m (\Omega)$ is the Marcinkiewicz space (see Definition 3.8 in \cite{BC}), which consists of all measurable functions $f$ on $\Omega$ with the following property: there exists a constant $\gamma$ such that
\begin{equation}\label{weak space}
|\{x\in \Omega: |f|>\lambda \}| \le \frac {\gamma}{ \lambda ^m}, \ \ \forall \lambda >0.
\end{equation}
The norm of $f\in L_{weak} ^m (\Omega)$ is defined by
$$
\|f\|_{L_{weak}^m (\Omega)} ^m =\inf \{\gamma>0, (\ref{weak space}) \mbox { hlods}\}.
$$
A H\"older inequality holds true for $f\in L_{weak}^m (\Omega)$, $m>1$: there exists $B=B( \|f\| _{L_{weak} ^m (\Omega)},m)>0$ such that for every measurable subset $E\subset \Omega$,
\begin{equation}\label{Holder}
\int_E |f| dx \le B |E| ^{1-\frac 1 m }.
\end{equation}

For $m>(2^*)'$ we consider weak solutions of (\ref{BVP}).

\begin{definition}
A function $u\in W_0 ^{1,2} (\Omega)$ is a solution to (\ref{BVP}) if for any $v\in W_0^{1,2} (\Omega)$,
\begin{equation}\label{weak solution}
\int_\Omega a(x,u) Du Dv dx =\int_\Omega f v dx.
\end{equation}
\end{definition}

We note that, because of assumption (\ref{conditon for a(x,s)}), the differential operator
$$
-\mbox {div} (a(x,v)Dv)
$$
is not coercive on $W_0^{1,2} (\Omega)$, even if it is well defined between $W_0^{1,2} (\Omega)$ and its dual, see the last chapter in the monograph \cite{BC}. We mention that, for $f\in L^m(\Omega)$, Boccardo, Dall'Aglio
and Orsina derived in \cite{BDO} some existence and regularity results of (\ref{BVP}), see also the last chapter in \cite{BC}.

In the first part of this section, we consider degenerate elliptic equations with the right hand side belongs to Marcinkiewicz space $L_{weak} ^m (\Omega)$ with $m>(2^*)'$ by using Lemma \ref{Stampacchia Lemma,
generalization}. This method seems to be simpler than the classical ones.

\begin{theorem}\label{theorem 1}
Let $f\in L_{weak} ^m (\Omega)$, $m>(2^*)'$, and $u\in W_0^{1,2} (\Omega)$ be a solution to (\ref{BVP}). Then

(i) $m>\frac n 2 \Rightarrow \exists L=L(n,\alpha, \|f\| _{L_{weak} ^m (\Omega)},m, |\Omega|)>0$ such that $|u| \le 2L$, a.e. $\Omega$;

(ii) $m =\frac n 2 \Rightarrow \exists \lambda =\lambda (n,\alpha, \|f\| _{L_{weak} ^m (\Omega)},m) >0 \mbox { such that }  e^{\lambda |u| ^{1-\theta }} \in L^1(\Omega) $;

(iii) $ (2^*)'< m<\frac n 2 \Rightarrow u\in L_{weak} ^{m^{**} (1-\theta)} (\Omega)$, $m^{**} =\frac {Nm}{N-2m}$.
\end{theorem}

If we weaken the summability hypotheses on $f$, then it is possible to give a meaning to solution for problem (\ref{BVP}), using the concept of entropy solutions
which has been introduced in \cite{BBGGPV}.

\begin{definition}\label{entropy solution}
Let $f\in L_{weak} ^m (\Omega)$, $1<m<(2^*)'$. A measurable function $u$ is an entropy solution of (\ref{BVP}) if $T_\ell (u)$ belongs to $W_0^{1,2} (\Omega)$
for every $\ell > 0$ and if
\begin{equation}\label{definition 2}
\int_\Omega a(x,u) Du DT_\ell (u-\varphi) dx \le \int_\Omega f T_ \ell (u-\varphi) dx ,
\end{equation}
for every $\ell >0$ and for every $\varphi \in W_0^{1,2} (\Omega) \cap L^\infty (\Omega)$.
\end{definition}

In (\ref{definition 2}),
$$
T_\ell (v) =\max \{-v,\min \{v,\ell \}\}
$$
is the truncation of $v$ at level $\ell >0$.

We have the following result for entropy solutions.

\begin{theorem}\label{theorem 2}
Let $u$ be an entropy solution of (\ref{BVP}) and $f\in L_{weak} ^m (\Omega)$, $1<m<(2^*)'$. Then $u\in L_{weak} ^{m^{**} (1-\theta)} (\Omega)$.
\end{theorem}

\vspace{3mm}

We now use Lemma \ref{Stampacchia Lemma, generalization} to prove Theorems \ref{theorem 1} and \ref{theorem 2}.

\vspace{3mm}

\noindent {\it Proof of Theorem \ref{theorem 1}.} Denote $G_k(u) =u-T_k(u)$.
For $h>k >0$, we take
$$
v =T_{h-k} (G_k (u))
$$
as a test function in (\ref{weak solution}) and we have, by (\ref{conditon for a(x,s)}), that
\begin{equation}\label{the first inequality}
\alpha \int_{B_{k,h}} \frac {|Du|^2}{ (1+|u|)^\theta} dx \le \int _{A_k} fT_{h-k} (G_k(u))dx,
\end{equation}
where
$$
B_{k,h} =\{x\in \Omega: k<|u| \le h\}, \ \ A_k =\{x\in \Omega: |u|>k\}.
$$
We now estimate both sides of (\ref{the first inequality}). The left hand side can be estimated as
\begin{equation}\label{the left hand side}
\begin{array}{llll}
&\displaystyle \alpha \int_{B_{k,h}} \frac {|Du| ^2}{ (1+|u|) ^\theta} dx \\
\ge &\displaystyle \frac {\alpha}{ (1+h) ^\theta} \int _{B_{k,h}} |Du|^2dx \\
= &\displaystyle \frac {\alpha}{ (1+h) ^\theta} \int _{\Omega} |DT_{h-k } (G_k(u))|^2dx \\
\ge &\displaystyle \frac {\alpha {\cal S}^2}{(1+h) ^\theta} \left(\int_\Omega \left|T_{h-k} (G_k(u))\right|^{2^*} \right) ^{\frac 2 {2^*}} \\
\ge &\displaystyle \frac {\alpha {\cal S}^2}{(1+h) ^\theta} \left(\int_{A_h} \left|T_{h-k} (G_k(u))\right|^{2^*} \right) ^{\frac 2 {2^*}} \\
= &\displaystyle \frac {\alpha {\cal S}^2}{(1+h) ^\theta} (h-k) ^2 |A_h| ^{\frac 2 {2^*}},
\end{array}
\end{equation}
where we have used the Sobolev inequality
$$
u\in W_0^{1,p} (\Omega) \Longrightarrow {\cal S} \|u\| _{L^{p^*}(\Omega)} \le \|Du\|_{L^p(\Omega)}, \ \ p\ge 1, \ {\cal S} ={\cal S} (n,p).
$$
Using H\"older inequality (\ref{Holder}), the right hand side can be estimated as
\begin{equation}\label{the right hand side}
\int_{A_k}  fT_{h-k} (G_k(u)) dx \le (h-k) \int_{A_k } |f| dx \le (h-k) B  |A_k| ^{\frac 1 {m'}},
\end{equation}
where $B$ is a constant depending on $\|f\|_{L_{weak}^m (\Omega)}$ and $m$. (\ref{the first inequality}) together with (\ref{the left hand side}) and (\ref{the right hand side}) implies, for $h\ge 1$,
$$
|A_h| \le \frac {c_6  (1+h) ^{\frac {\theta 2^*}{2}}}{ (h-k) ^{\frac {2^*}{2}}} |A_k| ^{\frac {2^*}{2m'}}\le \frac {c_7 h^{\frac {\theta 2^*}{2}}}{ (h-k) ^{\frac {2^*}{2}}} |A_k| ^{\frac {2^*}{2m'}},
$$
where $c_6$ is a constant depending on $n,\alpha, \|f\| _{L_{weak} ^m (\Omega)},m$ and $c_7 =c_6 2 ^{\frac {\theta 2^*}{2}}$. Thus (\ref{stampacchia_assumption 22}) holds true with
$$
\varphi (k) =|A_k|, \ c=c_7, \ \alpha =\frac {2^*}{2},  \ \beta =\frac {2^*}{2m'} \ \mbox { and } \ k_0=1.
$$
We use Lemma \ref{Stampacchia Lemma, generalization} and we have:

(i) If $m>\frac n 2$, then $\beta>1$. We use Lemma \ref{Stampacchia Lemma, generalization} (i) and we have $|A_{2L}| =0$ for some constant $L$ depending on $n,\alpha, \|f\| _{L_{weak} ^m (\Omega)},m$ and $|\Omega|$, from
which we derive $|u|\le 2L$ a.e. $\Omega$;

(ii)  If $m=\frac n 2$, then $\beta=1$. We use Lemma \ref{Stampacchia Lemma, generalization} (ii) and we derive that there exists a constant $\tau$ depending on $n,\alpha, \|f\| _{L_{weak} ^m (\Omega)}, m$ and $\theta$,
such that
$$
|\{|u|>k\}| \le |\{ |u|>1\}| e ^{1-\left(\frac {k-1}{\tau} \right) ^{1-\theta}} \le |\Omega| e ^{1-\left(\frac {k-1}{\tau} \right) ^{1-\theta}}.
$$
We let $2^{2-\theta} \lambda =\tau ^{\theta -1}$ and $k\ge 2$ ($\Leftrightarrow k-1\ge \frac k 2$) and we have
$$
|\{|u|>k\}| \le |\Omega| e ^{1-2^{2-\theta}\lambda (k-1)^{1-\theta}} \le |\Omega| e ^{1-2^{2-\theta}\lambda \left(\frac k 2 \right) ^{1-\theta}}=c_7 e ^{-2\lambda k ^{1-\theta}}, \ k\ge 2,
$$
where $c_7= |\Omega| e $. Hence
$$
|\{e ^{\lambda |u| ^{1-   \theta}} >e ^{\lambda k ^{1-   \theta}}\}| = |\{ |u| >k\}| \le c _7 e ^{-2\lambda k ^{1-   \theta}}.
$$
Let $ \tilde k =e ^{\lambda k ^{1-   \theta}}$, then
$$
|\{e ^{\lambda |u| ^{1-   \theta}} >  \tilde  k \}| \le \frac {c_7}{  \tilde k ^2}, \ \mbox { for } \  \tilde k \ge e ^{\lambda 2 ^{1-\theta}}.
$$
Let $\big [e ^{\lambda 2 ^{1-\theta}}\big] =k_1$.  We now use Lemma 3.11 in \cite{BC} which states that the sufficient and necessary condition for $g\in L^r (\Omega)$, $r\ge 1$, is
$$
\sum_{k=0}^\infty k^{r-1} |\{|g| >k\}| <+\infty.
$$
We use the above lemma for $g= e ^{\lambda |u| ^{1-   \theta}}$ and $r=1$. Since
$$
\sum_{\tilde k=0} ^\infty |\{e ^{\lambda |u| ^{1-   \theta}} >\tilde   k \}|= \left (\sum_{\tilde  k=0} ^ {k_1}+ \sum _{\tilde k=k_1+1} ^\infty\right) \le (k_1+1) |\Omega| + c _7\sum _{\tilde k=k_1+1} ^\infty \frac 1 {
\tilde  k ^2} <+\infty,
$$
then $e ^{\lambda |u| ^{1-\theta}} \in L^1 (\Omega)$;

(iii) If $(2^*)'< m<\frac n 2$, then $\beta<1$. We use Lemma \ref{Stampacchia Lemma, generalization} (iii) and we have for all $k\ge 1$,
$$
|\{|u| >k\}| \le c_8 \left(\frac 1 k \right) ^{\frac {(1-\theta) \alpha}{1-\beta}} = c_8 \left(\frac 1 k \right) ^{m^{**} (1-\theta)} ,
$$
where $c_8$ is a constant depending on $n,\alpha, \|f\| _{L_{weak} ^m (\Omega)},m$ and $|\Omega|$. Thus
$$
|\{|u| >\lambda \}| \le \max\{ c_8, |\Omega|\} \left(\frac 1 \lambda \right) ^{m^{**} (1-\theta)} , \ \forall \lambda >0,
$$
that is $u\in L_{weak} ^{m^{**} (1-\theta)} (\Omega)$, as desired.     \qed

\vspace{4mm}

\noindent {\it Proof of Theorem \ref{theorem 2}.}  Let $h>k>0$. In (\ref{definition 2}) we take $\varphi =T_k(u) \in W_0^{1,2} (\Omega)\cap L^\infty (\Omega)$ and $\ell =h-k$.
By (\ref{conditon for a(x,s)}) we derive (\ref{the first inequality}). The result follows from the proof of Theorem \ref{theorem 1} (iii).    \qed

\vspace{3mm}

\noindent
{\bf Acknowledgments.}  The first was supported by NSFC (12071021), NSF of Hebei Province (A2019201120) and Hebei Provincial Key Project of Science and Technology (ZD2021307); the third author was supported by the Science
and Technology Research Youth Fund Project of Higher Education Institutions of Hebei Province (QN2020124).

\rm 
\end{document}